\documentclass[11pt]{article}

\usepackage{graphicx}        % standard LaTeX graphics tool
\makeindex             % used for the subject index

\usepackage{graphicx}

\def\be{\begin{equation}}
\def\ee{\end{equation}}
\newcommand{\kk}[2]{\frac{#1}{#2}}

\newcommand{\ff}[1]{{\bf  #1}}
\def\a{\alpha}
\def\b{\beta}

\def\e{\epsilon}
\def\lam{\lambda}
\def\ra{\rightarrow}
\def\={\approx}
\def\x{\ff{x}}
\def\vcode#1#2#3#4{\begin{figure}
\begin{center}
\begin{minipage}[c]{#1\textwidth}
{{\small #2 \hrule  \vspace{5pt}   %
{\it #3} \hrule }}
\end{minipage}
\caption{#4}
\end{center}   \end{figure}} %%
\begin{document}

\title{Eagle Strategy Using L\'evy Walk and Firefly Algorithms For Stochastic Optimization}

\author{Xin-She Yang\footnote{Corresponding Author}  \\
Department of Engineering,University of Cambridge, \\
Trumpinton Street, Cambridge CB2 1PZ, UK
% \email{xy227@cam.ac.uk}
\and
Suash Deb, Department of Computer Science \& Engineering, \\
C. V. Raman College of Engineering,
Bidyanagar, Mahura, Janla, \\
Bhubaneswar 752054, INDIA
% \email{suashdeb@gmail.com}
}

\date{}

\maketitle

%% Abstract

\abstract{Most global optimization problems are nonlinear and thus difficult to solve,
and they become even more challenging when uncertainties are present in objective
functions and constraints. This paper provides a new two-stage hybrid search method, called
Eagle Strategy, for stochastic optimization. This strategy intends to combine the
random search using L\'evy walk with the firefly algorithm in an iterative manner.
Numerical studies and results suggest that the proposed Eagle Strategy is very
efficient for stochastic optimization.  Finally practical
implications and potential topics for further research will be discussed. \\
}

%%
%% {\bf Key words:} Eagle strategy; firefly algorithm; metaheuristic; stochastic optimization.
%%

\noindent {\bf Citation detail:} 

X.-S. Yang and S. Deb, Eagle strategy using Levy walk and firefly algorithms for
stochastic optimization, in: {\it Nature Inspired Cooperative Strategies for Optimization
(NISCO 2010)} (Eds. J. R. Gonzalez et al., Studies in Computational Intelligence, Springer Berlin, {\bf 284}, 101-111 (2010).

%% Begin of Main Text %%

\section{Introduction}

To find the solutions to any optimization problems, we can use  either conventional
optimization algorithms such as the Hill-climbing and simplex method, or heuristic
methods such as genetic algorithms, or their proper combinations.
Modern metaheuristic algorithms are becoming powerful in solving
global optimization problems \cite{Bod,Deb,Gold,Ken2,Yang,Yang2,Yang3},
especially for the NP-hard problems such
as the travelling salesman problem. For example,
particle swarm optimization (PSO) was developed by Kennedy and
Eberhart in 1995 \cite{Ken,Ken2}, based on the swarm behaviour such
as fish and bird schooling in nature. It has now been applied
to find solutions for many optimization applications.
Another example is the Firefly Algorithm developed by the first author
\cite{Yang,Yang2} which has demonstrated promising superiority over
many other algorithms. The search strategies in these multi-agent
algorithms are controlled randomization and exploitation of the
best solutions. However, such randomization typically uses
a uniform distribution or Gaussian distribution. In fact, since
the development of PSO, quite a few algorithms have been developed
and they can outperform PSO in different ways \cite{Yang2,YangYao}.

On the other hand, there is always some uncertainty and noise associated with
all real-world optimization problems. Subsequently,
objective functions may have noise and constraints may also have
random noise. In this case, a standard optimization problem
becomes a stochastic optimization problem. Methods that work well
for standard optimization problems cannot directly be applied
to stochastic optimization; otherwise, the obtained
results are incorrect or even meaningless. Either the optimization
problems have to be reformulated properly or the optimization
algorithms should be modified accordingly, though in most cases
we have to do both \cite{Bental,Marti,Wallace}

In this paper, we intend to formulate a new metaheuristic
search method, called Eagle Stategy (ES), which combines
the L\'evy walk search with the Firefly Algorithm (FA).
We will provide the comparison study of the ES with PSO and
other relevant algorithms. We will first outline the basic
ideas of the Eagle Strategy, then outline the essence of
the firefly algorithm, and finally carry out the comparison
about the performance of these algorithms.

\section{Stochastic Multiobjective  Optimization}

An ordinary optimization problem, without noise or uncertainty,
can be written as
\be \min_{\x \in \Re^d} f_i(\x), \; (i=1,2,...,N)  \ee
\[ \textrm{subject to } \phi_j(\x)=0, \; (j=1,2,...,J), \]
\be \qquad \qquad \psi_k(\x) \le 0,  \;  (k=1,2,...,K), \ee
where $\x=(x_1,x_2, ..., x_d)^T$ is the vector of design
variables.

For stochastic optimization problems, the effect of uncertainty or noise on the design
variable $x_i$ can be described by a random variable $\xi_i$ with
a distribution $Q_i$ \cite{Marti,Wallace}. That is
\be x_i \mapsto \xi_i(x_i), \ee
and
\be \xi_i \sim Q_i. \ee
The most widely used distribution
is the Gaussian or normal distribution $N(x_i, \sigma_i)$ with
a mean $x_i$ and a known standard deviation $\sigma_i$. Consequently,
the objective functions $f_i$ become random variables $f_i(\x,\xi)$.

Now we have to reformulate the optimization problem as the minimization of the mean of the
objective function $f_i(\x)$ or $\mu_{f_i}$
\be \min_{\x \in \Re^d} \{ \mu_{f_1}, ..., \mu_{f_N} \}. \ee
Here $\mu_{f_i}=E(f_i)$ is the mean or expectation of $f_i(\xi(\x))$ where $i=1,2,...,N$.
More generally, we can also include their uncertainties, which leads to the
minimization of
\be \min_{\x \in \Re^d} \{ \mu_{f_1} + \lam \sigma_1, ..., \mu_{f_N} + \lam \sigma_N \}, \ee
where $\lam \ge 0$ is a constant. In addition, the constraints with uncertainty
should be modified accordingly.

In order to estimate $\mu_{f_i}$, we have to use some sampling techniques such as the Monte Carlo
method. Once we have randomly drawn the samples, we have
\be \mu_{f_i} \= \kk{1}{N_i} \sum_{p=1}^{N_i} f_i(\x, \xi^{(p)}), \ee
where $N_i$ is the number of samples.

\section{Eagle Strategy}

The foraging behaviour of eagles such as golden eagles or {\it Aquila Chrysaetos}
is inspiring. An eagle forages in its own territory by flying freely in a random
manner much like the L\'evy flights. Once the prey is sighted, the eagle will
change its search strategy to an intensive chasing tactics so as to catch the
prey as efficiently as possible. There are two important components to an eagle's
hunting strategy: random search by L\'evy flight (or walk)
and intensive chase by locking its  aim on the target.

Furthermore, various studies have shown that
flight behaviour of many animals and insects has demonstrated
the typical characteristics of L\'evy flights
\cite{Brown,Reynolds,Pav,Pav2}. A recent study by Reynolds and Frye shows that
fruit flies or {\it Drosophila melanogaster}, explore their landscape using a series
of straight flight paths punctuated by a sudden $90^{0}$ turn, leading to
a L\'evy-flight-style intermittent scale-free search pattern.
Studies on human behaviour such as the Ju/'hoansi hunter-gatherer foraging
patterns also show the typical feature of L\'evy flights.
Even light can be related to L\'evy flights \cite{Barth}.
Subsequently, such behaviour has been applied to
optimization and optimal search, and preliminary results show its
promising capability \cite{Pav,Reynolds,Shles,Shles2}.

\subsection{Eagle Strategy}

Now let us idealize the two-stage strategy of an eagle's foraging behaviour.
Firstly, we assume that an eagle will perform the L\'evy walk in the whole domain.
Once it finds a prey it changes to a chase strategy. Secondly, the chase strategy
can be considered as an intensive local search using any optimization technique such
as the steepest descent method, or the downhill simplex or Nelder-Mead method \cite{Nelder}.
Obviously, we can also use any efficient metaheuristic algorithms
such as the particle swarm optimization (PSO)
and the Firefly Algorithm (FA) to do concentrated local search. The pseudo code of the proposed
eagle strategy is outlined in Fig. \ref{fa-fig-100}.

The size of the hypersphere depends on the landscape of the objective functions.
If the objective functions are unimodal, then the size of the hypersphere can be about the
same size of the domain. The global optimum can in principle be found from any initial guess.
If the objective are multimodal, then the size of the hypersphere should be the typical
size of the local modes. In reality, we do not know much about the landscape of the objective
functions before we do the optimization, and we can either start from a larger domain
and shrink it down or use a smaller size and then gradually expand it.

On the surface, the eagle strategy has some similarity with the random-restart hill climbing
method, but there are two important differences. Firstly, ES is a two-stage strategy rather
than a simple iterative method, and thus ES intends to combine a good randomization (diversification) technique
of global search with an intensive and efficient local search method. Secondly, ES uses L\'evy
walk rather than simple randomization, which means that the global search space can be
explored more efficiently. In fact, studies show that L\'evy walk is far more efficient
than simple random-walk exploration.

\vcode{0.7}{{\sf Eagle Strategy}}{
\indent Objective functions $f_1(\x),...,f_N(\x)$ \\
\indent \quad Initial guess $\x^{t=0}$  \\
\indent \quad {\bf while} ($||\x^{t+1}-\x^t||>$ tolerance) \\
\indent \quad Random search by performing L\'evy walk \\
\indent \quad Evaluate the objective functions \\
\indent \quad Intensive local search with a hypersphere \\
\indent \quad \qquad via Nelder-Mead or the Firefly Algorithm \\
\indent \qquad  {\bf if} (a better solution is found) \\
\indent \qquad \quad Update the current best \\
\indent \qquad  {\bf end if} \\
\indent \qquad  Update $t=t+1$ \\
\indent \qquad  Calculate means and standard deviations \\
\indent \quad {\bf end while } \\
\indent Postprocess results and visualization }{Pseudo code of the Eagle Strategy (ES).
\label{fa-fig-100} }

The L\'evy walk has a random step length
being drawn from a L\'evy distribution
\be \textrm{L\'evy} \sim u = t^{-\lam}, \quad (1 < \lam \le 3), \ee
which has an infinite variance with an infinite mean.
Here the steps of the eagle motion is essentially
a random walk process with a power-law step-length distribution with a heavy tail.
The special case $\lam=3$ corresponds to Brownian motion, while
$\lam=1$ has a characteristics of stochastic tunneling, which may be more efficient
in avoiding being trapped in local optima.

For the local search, we can use any efficient optimization algorithm such as the
downhill simplex (Nelder-Mead) or metaheuristic algorithms such as PSO and
the firefly algorithm. In this paper, we used
the firefly algorithm to do the local search, since the firefly algorithm was designed to
solve multimodal global optimization problems \cite{Yang2}.

\subsection{Firefly Algorithm}

We now briefly outline the main components of the Firefly Algorithm developed
by the first author \cite{Yang}, inspired by the flash pattern and
characteristics of fireflies.  For simplicity in
describing the algorithm, we now use the following
three idealized rules: 1) all fireflies are unisex so that one firefly will be attracted to other fireflies
regardless of their sex; 2) Attractiveness is proportional to their brightness,
thus for any two flashing fireflies, the less brighter one will move towards the brighter
one. The attractiveness is proportional to the brightness and they both
decrease as their distance increases. If there is
no brighter one than a particular firefly, it will move randomly;
3) The brightness of a firefly is affected or  determined by the landscape of the
objective function. For a  maximization problem, the brightness can simply be proportional
to the value of the objective functions.

\vcode{0.77}{{\sf Firefly Algorithm}} {
\indent Objective function $f_p(\x), \; \x=(x_1, ..., x_d)^T$ \\
\indent Initial population of fireflies $\x_i \; (i=1,...,n)$ \\
\indent Light intensity $I_i$ at $\x_i$ is determined by $f_p(\x_i)$  \\
\indent Define light absorption coefficient $\gamma$ \\
\indent {\bf while} ($t<$MaxGeneration) \\
\indent \quad {\bf for} $i=1:n$ all $n$ fireflies \\
\indent \qquad {\bf for} $j=1:i$ all $n$ fireflies \\
\indent \qquad {\bf if} ($I_j>I_i$) \\
\indent \qquad  \qquad Move firefly $i$ towards $j$ ($d$-dimension) \\
\indent \qquad  {\bf end if} \\
\indent \qquad  Vary $\beta$ via $\exp[-\gamma r]$ \\
\indent \qquad  Evaluate new solutions and update \\
\indent \qquad {\bf end for }$j$ \\
\indent \quad {\bf end for }$i$ \\
\indent Rank the fireflies and find the current best  \\
\indent {\bf end while} \\
\indent  Postprocess results and visualization }{Pseudo code of the firefly algorithm (FA).
\label{fa-fig-200} }

In the firefly algorithm, there are two important issues: the variation of
light intensity and formulation of the attractiveness.
For simplicity, we can always assume that the attractiveness
of a firefly is determined by its brightness which in turn is associated with
the encoded objective function.

In the simplest case for maximum optimization problems,
the brightness $I$ of a firefly at a particular location $\x$ can be chosen
as $I(\x) \propto f(\x)$. However, the attractiveness $\b$ is relative, it should be
seen in the eyes of the beholder or judged by the other fireflies. Thus, it will
vary with the distance $r_{ij}$ between firefly $i$ and firefly $j$. In addition,
light intensity decreases with the distance from its source, and light is also
absorbed in the media,
so we should allow the attractiveness to vary with the degree of absorption.
In the simplest form, the light intensity $I(r)$ varies according to the inverse
square law $I(r)=\kk{I_s}{r^2}$ where $I_s$ is the intensity at the source.
For a given medium with a fixed light absorption coefficient $\gamma$,
the light intensity $I$ varies with the distance $r$. That is
\be I=I_0 e^{-\gamma r}, \ee
where $I_0$ is the original light intensity.

As a firefly's attractiveness is proportional to
the light intensity seen by adjacent fireflies, we can now define
the attractiveness $\b$ of a firefly by
\be \b = \b_0 e^{-\gamma r^2}, \label{att-equ-100} \ee
where $\b_0$ is the attractiveness at $r=0$.

The distance between any two fireflies $i$ and $j$ at $\x_i$ and $\x_j$, respectively, is
the Cartesian distance
\be r_{ij}=||\x_i-\x_j|| =\sqrt{\sum_{k=1}^d (x_{i,k} - x_{j,k})^2}, \ee
where $x_{i,k}$ is the $k$th component of the spatial coordinate $\x_i$ of $i$th
firefly. In the 2-D case, we have \be r_{ij}=\sqrt{(x_i-x_j)^2+(y_i-y_j)^2}.\ee

The movement of a firefly $i$ is attracted to another more attractive (brighter)
firefly $j$ is determined by
\be \x_i =\x_i + \b_0 e^{-\gamma r^2_{ij}} (\x_j-\x_i) + \a \; ({\rm rand}-\kk{1}{2}), \ee
where the second term is due to the attraction.  The third term
is randomization with a control parameter $\a$, which makes the exploration of the search
space more efficient.

We have tried to use different values of the parameters $\a, \b_0, \gamma$ \cite{Yang,Yang2},
after some simulations, we concluded that we can
use $\b_0=1$, $\a \in [0,1]$, $\gamma=1$, and $\lam=1$ for most applications.
In addition, if the scales vary significantly in different dimensions such as $-10^5$ to $10^5$ in
one dimension while, say, $-0.001$ to $0.01$ along the other, it is a good idea to
replace $\a$ by $\a S_k$ where the scaling parameters $S_k (k=1,...,d)$ in
the $d$ dimensions should be determined by the actual scales of the problem of interest.

There are two important limiting cases when $\gamma \ra 0$ and $\gamma \ra \infty$.
For $\gamma \ra 0$, the attractiveness is constant $\b=\b_0$ and the length scale $\Gamma=1/\sqrt{\gamma} \ra \infty$,
this is equivalent to say that the light intensity does not decrease in an idealized sky.
Thus, a flashing firefly can be seen anywhere in the domain. Thus, a single (usually global)
optimum can easily be reached. This corresponds to a special case of particle
swarm optimization (PSO) discussed earlier. Subsequently, the efficiency of this special case
could be about the same as that of PSO.

On the other hand, the limiting case $\gamma \ra \infty$ leads to
$\Gamma \ra 0$ and  $\b(r) \ra \delta(r)$ (the Dirac delta function), which means that
the attractiveness is almost zero in the sight of other fireflies or the fireflies are short-sighted. This is equivalent
to the case where the fireflies roam in a very foggy region randomly. No other fireflies can be
seen, and each firefly roams in a completely random way. Therefore, this corresponds
to the completely random search method.
As the firefly algorithm is usually in somewhere between these two extremes, it is possible
to adjust the parameter $\gamma$ and $\alpha$ so that it can outperform both the random search and PSO.

\section{Simulations and Comparison}

\subsection{Validation}

\begin{figure}
 \centerline{\includegraphics[width=5in,height=3.5in]{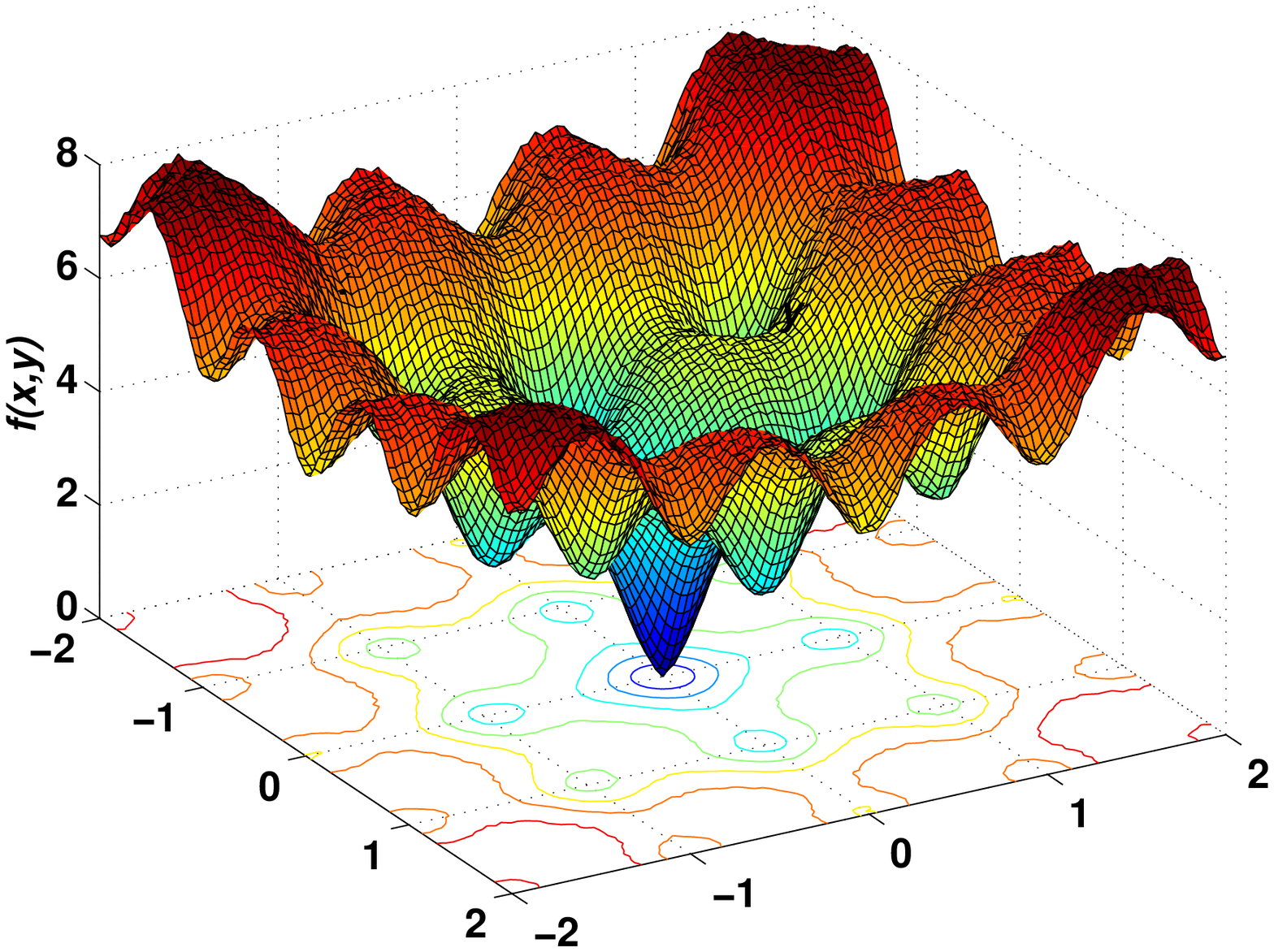} }
\vspace{-10pt}
\caption{Ackley's 2D function with Gaussian noise. \label{yangfa-200} }
\end{figure}

In order to validate the proposed algorithm, we have implemented it in
Matlab.  In our simulations,  the values
of the parameters are $\a=0.2$, $\gamma=1$, $\lam=1$, and $\b_0=1$.
\begin{figure}
 \centerline{\includegraphics[width=5in,height=3.5in]{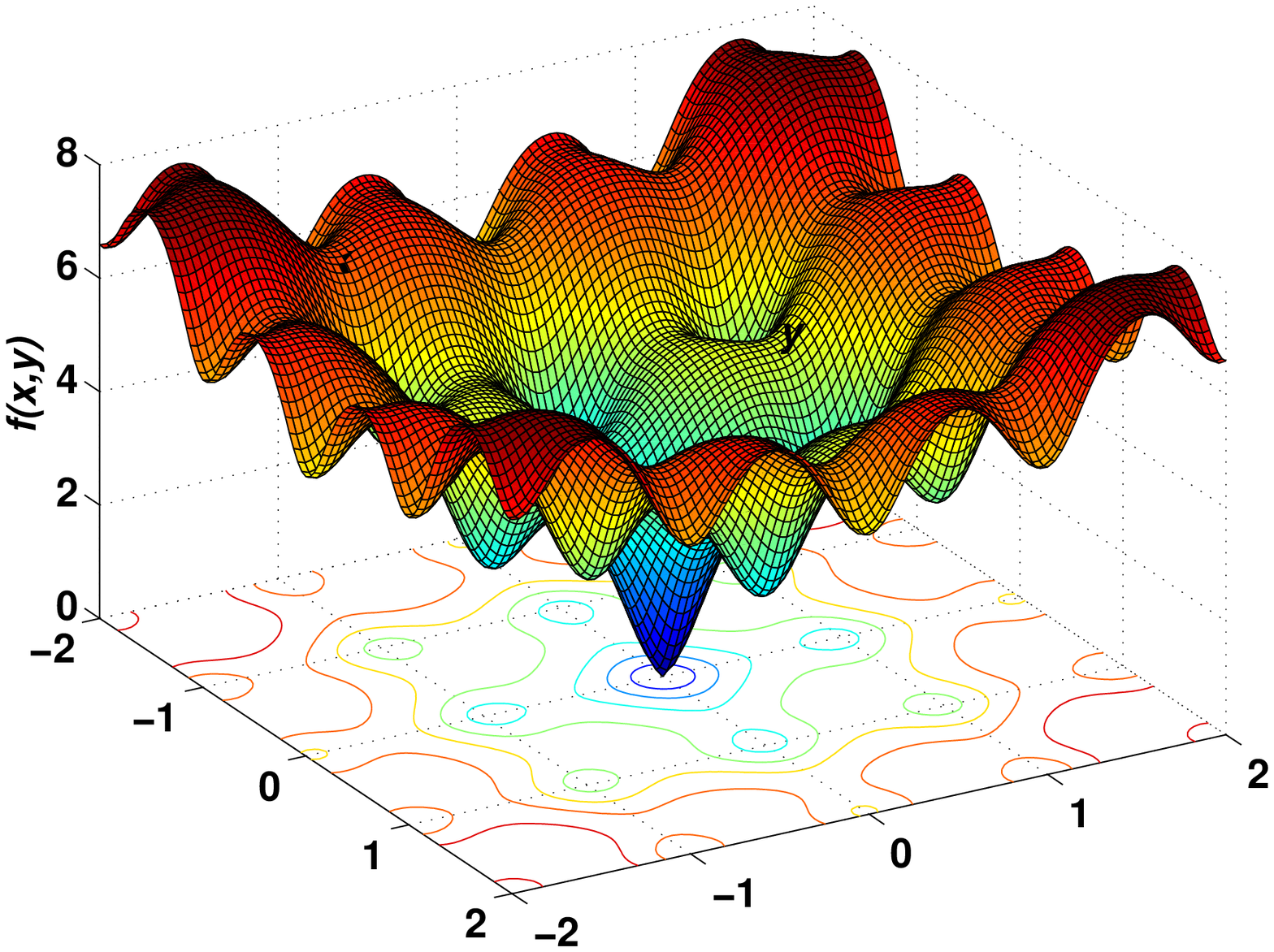} }
\vspace{-10pt}
\caption{Ackley's function for two independent variables
with a global minimum $f_*=0$ at $(0,0)$. \label{yangfa-100} }
\end{figure}
As an example, we now use the ES to find the
global optimum of the Ackley function
\be f(\x)=-20 \exp[-\kk{1}{5} \sqrt{\kk{1}{d} \sum_{i=1}^d x_i^2} ]
 -\exp[\kk{1}{d} \sum_{i=1}^d \cos(2 \pi x_i)]+20+e,  \ee
where $(d=1,2,...)$ \cite{Ackley}. The global minimum $f_*=0$  occurs at $(0,0,...,0)$
in the domain of $-32.768 \le x_i \le 32.768$ where $i=1,2,...,d$.
The landscape of the 2D Ackley function is shown in Fig. \ref{yangfa-100},
while the landscape of this function with $2.5\%$ noise
is shown in Fig. \ref{yangfa-200}

\begin{figure}
 \centerline{\includegraphics[height=2.5in,width=3in]{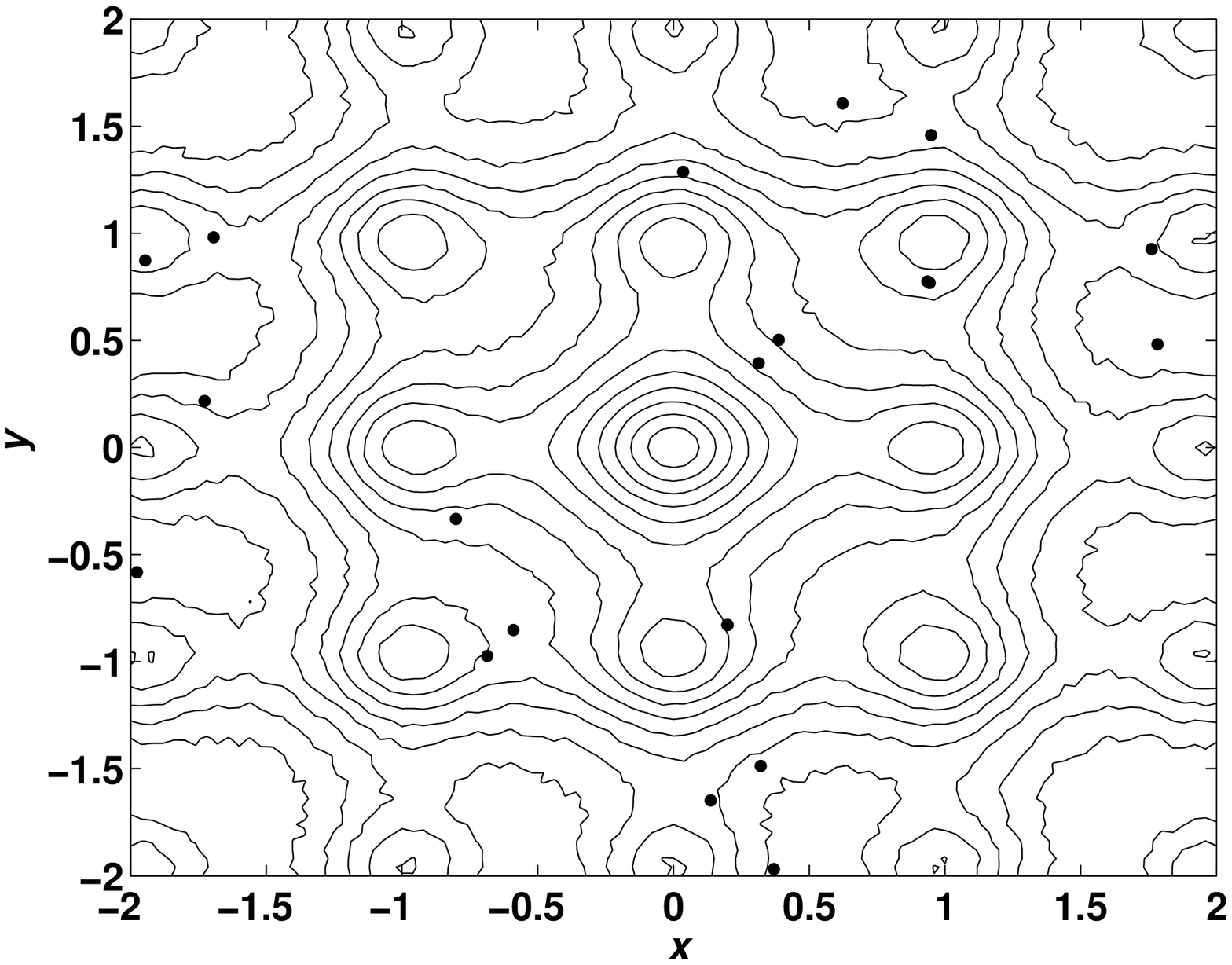}
 \includegraphics[height=2.5in,width=3in]{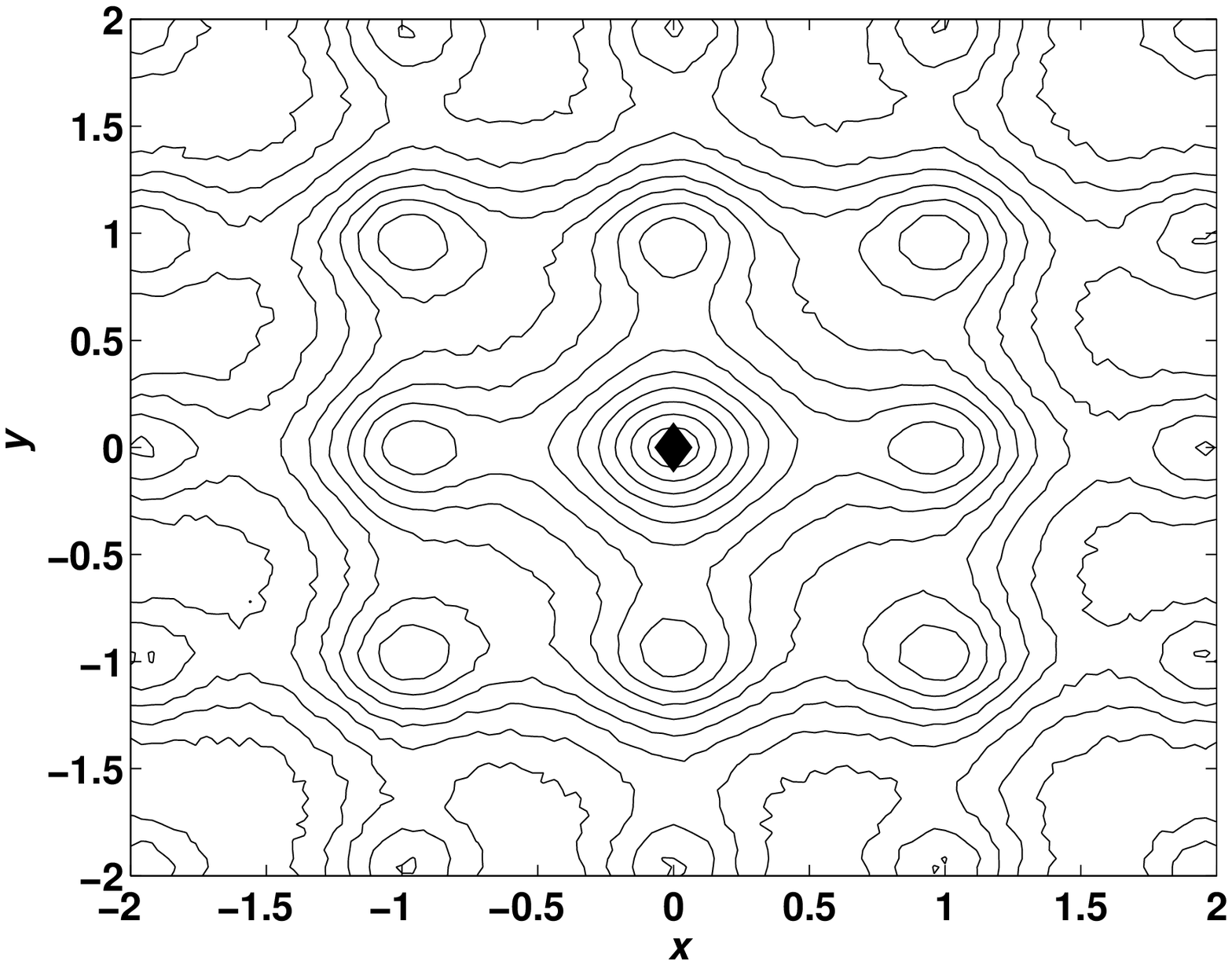}}
\vspace{-5pt}
\caption{The initial locations of the 20 fireflies (left) and their
locations after $15$ iterations (right). We have used $\gamma=1$. \label{yangfa-300} }
\end{figure}

The global minimum in 2D for a given noise level of $2.5\%$
can be found after about 300 function evaluations (for 20 fireflies after 15
iterations, see Fig. \ref{yangfa-300}).

\subsection{Comparison of ES with PSO}

Various studies show that PSO algorithms can
outperform genetic algorithms (GA) \cite{Gold} and other conventional algorithms for
solving many optimization problems. This is
partially due to that fact that the broadcasting ability of the current
best estimates gives better and quicker convergence towards the
optimality. A general framework for evaluating statistical performance of evolutionary algorithms
has been discussed in detail by Shilane et al. \cite{Shilane}.

Now we will compare the Eagle Strategy with PSO
for various standard test functions.
After implementing these algorithms using
Matlab, we have carried out extensive simulations and each algorithm has been
run at least 100 times so as to carry out meaningful statistical analysis.
The algorithms stop when the variations of function values are less than
a given tolerance $\e \le 10^{-5}$. The results are
summarized in the following table (see Table 1) where the global optima
are reached. The numbers are in the
format: average number of evaluations (success rate), so
$12.7 \pm 1.15 (100)$ means that the average number (mean) of function
evaluations is $12.7 \times 10^3=12700$ with a standard deviation of $1.15 \times 10^3=1150$.
The success rate of finding the global optima for this algorithm is $100\%$.
Here we have used the following abbreviations:  MWZ for Michalewicz's function
with $d=16$, RBK for Rosenbrock with $d=16$, De Jong for De Jong's sphere function
with $d=256$,
Schwefel for Schwefel with $d=128$, Ackley for Ackley's function with $d=128$,
and Shubert for Shubert's function with 18 minima. In addition, all these
test functions have a $2.5\%$ of Gaussian noise, or $\sigma=0.025$.
In addition, we have used the population size $n=20$ in all our simulations.
\begin{table}[ht]
\caption{Comparison of algorithm performance}
\centering
\begin{tabular}{ccccc}
\hline \hline
  &  PSO ($\times 10^3$)  & ES ($\times 10^3$)\\
\hline
 Easom  & $185.9 \pm  3.1 (97)$ & $12.7 \pm 1.15 (100 )$ \\
 MWZ  & $346.1 \pm 8.0 (98)$  & $36.1 \pm 3.5 (100)$ \\
 Rosenbrock   & $1637 \pm 79 (98)$ & $75 \pm 6.4 (100) $ \\
 De Jong & $852 \pm 16 (100)$ & $70.7 \pm 7.3 (100)$\\
 Schwefel & $726.1 \pm 25 (97)$ & $99 \pm 6.7 (100)$ \\
 Ackley  & $1170 \pm 19 (92)$ & $54 \pm 5.2 (100)$ \\
 Rastrigin & $3973 \pm 64 (90)$ & $151 \pm 14 (100)$ \\
 Easom & $863.7 \pm 55 (90)$ & $76 \pm 11 (100)$ \\
 Griewank  & $2798 \pm 63 (92)$ & $134 \pm 9.1 (100)$ \\
 Shubert  & $1197 \pm 56 (92)$ & $32 \pm 2.5 (100)$ \\

\hline
\end{tabular}
\end{table}

We can see that the ES is noticeably
more efficient in finding the global optima
with the success rates of  $100\%$. Each function evaluation is virtually instantaneous
on a modern personal computer. For example, the computing time for 10,000 evaluations
on a 3GHz desktop is about 5 seconds. Even with graphics for displaying the
locations of the particles and fireflies, it usually takes less than a few minutes.
Furthermore, we have used various values of the population size $n$ or the number of
fireflies. We found that for most problems $n=15$ to $50$ would be sufficient.
For tougher problems, larger $n$ such as $n=100$ or $250$ can be used, though excessively
large $n$ should not be used unless there is no better alternative,
as it is more computationally extensive.

\section{Conclusions}

By combining L\'evy walk with the firefly algorithm, we have successfully formulated
a hybrid optimization algorithm, called Eagle Strategy, for stochastic optimization.
After briefly outlining the basic procedure and its
similarities and differences with particle swarm optimization,
we then implemented and compared these algorithms.
Our simulation results for finding the global optima
of various test functions suggest that ES can significantly outperform
the PSO in terms of both efficiency and success rate.
This implies that ES is potentially more powerful in
solving NP-hard problems.

However, we have not carried out sensitivity studies of the algorithm-dependent
parameters such as the exponent $\lam$ in L\'evy distribution and the light absorption
coefficient $\gamma$,
which may be fine-tuned to a specific problem. This can form an important
research topic for further research. Furthermore, other local search algorithms
such as the Newton-Raphson method, sequential quadratic programming and Nelder-Mead
algorithms can be used to replace the firefly algorithm, and a comparison study
should be carried out to evaluate their performance. It may also show interesting
results if the level of uncertainty varies and it can be expected that the higher level
of noise will make it more difficult to reach optimal solutions.

As other important further studies, we can also focus on the applications of
this hybrid algorithm on the NP-hard traveling salesman problem. In addition, many engineering
design problems typically have to deal with intrinsic inhomogeneous materials properties and
such uncertainty may often affect the design choice in practice. The application of
the proposed hybrid algorithm in engineering design optimization may prove fruitful.

%% End of text %%

\end{document}